\documentclass[conference,final]{IEEEtran}
\IEEEoverridecommandlockouts

\usepackage{geometry}
\usepackage{amsmath,amssymb, theorem}
\usepackage{tikz}
\usepackage{epstopdf}
\usepackage{verbatim}
\usepackage{upgreek}
\usetikzlibrary{automata,positioning}
\pgfdeclarelayer{bg}    
\pgfsetlayers{bg,main}  
\pgfdeclarelayer{background}
\pgfdeclarelayer{foreground}
\pgfsetlayers{background,main,foreground}
\usepackage{ifthen}
\usepackage{psfrag}
\usepackage{cite}
\usepackage{graphicx}
\usepackage{graphics}
\usepackage{epstopdf}
\usepackage{bbm}
\usepackage{float}
\usepackage{framed}
\usepackage{textcomp}
\usepackage{mathabx}
\usepackage{subfigure}
\usepackage[font={small}]{caption}
\usepackage{accents}
\restylefloat{figure}
\usepackage{color,xspace}
\usepackage{tikz}
\usetikzlibrary{decorations.pathreplacing}
\usepackage{upgreek}

\usepackage{lipsum}

\usepackage{algorithm,algpseudocode}

\usetikzlibrary{automata,positioning}

\usepackage{pstool}
\usepackage{url}

\newcommand{\Kihx}{K_i^{\{0x\}}}
\newcommand{\Kihh}{K_i^{\{00\}}}
\newcommand{\Kihi}{K_i^{\{01\}}}
\newcommand{\bE}{\mathbb{E}}

\newcommand{\G}{G} 

\newcommand{\Nin}{\mathcal{N}^{\text{in}}}

\newcommand{\V}{V} 

\newcommand{\dmax}{d_{\max}}
\newcommand{\mor}{\operatorname{mor}}
\newcommand{\NP}{\mathcal{NP}}

\newcommand{\xhat}{\hat{x}}
\newcommand{\Xhatitt}{\hat{x}_i(t|t)}
\newcommand{\Xhatitmotmo}{\hat{x}_i (t-1|t-1)}
\newcommand{\Xhatitpot}{\hat{x}_i (t+1|t)}

\newcommand{\psiiw}{\psi_i^{\{w\}}}

\newcommand{\Hst}{\mathcal{H}}
\newcommand{\Obs}{\mathcal{O}}
\newcommand{\C}{\mathcal{C}}
\newcommand{\I}{\mathcal{I}}

\newcommand{\Emor}{\mor(E)}

\newtheorem{defn}{\textbf{Definition \newline}}

\newtheorem{theorem}{\textbf{Theorem \newline}}

\newcommand{\margin}[1]{\marginpar{\color{red}\tiny\ttfamily#1}}
\newcommand{\marginnw}[1]{\marginpar{\color{blue}\tiny\ttfamily#1}}

\newcommand{\oprocendsymbol}{\hbox{$\bullet$}}
\newcommand{\oprocend}{\relax\ifmmode\else\unskip\hfill\fi\oprocendsymbol}

\newtheorem{remark}{\textbf{Remark}}{}
{}

\begin{document}
\title{\huge Inference, Prediction and Control of Networked Epidemics}
\author{Nicholas J. Watkins, Cameron Nowzari, and George J. Pappas\thanks{N.J. Watkins and G.J. Pappas are with the Department of Electrical and Systems Engineering, University of Pennsylvania, Pennsylvania, PA 19104, USA, {\tt\small \{nwatk,pappasg\}@upenn.edu}; C. Nowzari is with the Department of Electrical and Computer Engineering, George Mason University, Fairfax,VA 22030, USA {\tt \small cnowzari@gmu.edu. } }}
\maketitle
\begin{abstract}  
We develop a feedback control method for networked epidemic spreading processes.  In contrast to most prior works which consider mean field, open-loop control schemes, the present work develops a novel framework for feedback control of epidemic processes which leverages incomplete observations of the stochastic epidemic process in order to control the exact dynamics of the epidemic outbreak. 
We develop an observation model for the epidemic process, and demonstrate that if the set of observed nodes is sufficiently well structured, then the random variables which denote the process' infections are conditionally independent given the observations.  We then leverage the attained conditional independence property to construct tractable mechanisms for the inference and prediction of the process state, avoiding the need to use mean field approximations or combinatorial representations.  We conclude by formulating a one-step lookahead controller for the discrete-time Susceptible-Infected-Susceptible (SIS) epidemic process which leverages the developed Bayesian inference and prediction mechanisms, and causes the epidemic to die out at a chosen rate.
\end{abstract}

\section{Introduction}

Networked epidemic spreading processes have been a subject of intense focus in the controls community for the past several years.  With envisioned applications ranging across disparate fields from the efficient mitigation of biological epidemics \cite{keeling2005networks}, to the efficient use of limited resources in ad-hoc communication networks \cite{pelusi2006opportunistic}, to analyzing the effectiveness of viral advertisement campaigns \cite{leskovec2007dynamics}, and the spread of innovations \cite{rogers2010diffusion}, there is little wonder as to why the control of such processes have garnered significant attention.  While a recent review of the state of the field can be found in \cite{nowzari2016analysis}, we now review the works most relevant to that which is presented here for purposes of completeness.

The majority of works to date considering networked epidemic processes study continuous-time epidemic processes under \emph{mean field} approximation, as direct analysis of the stochastic process is widely considered to be intractable.  Many efforts in this area use elements from spectral graph theory to pose convexification schemes for mean field resource allocation problems, the solutions of which guarantee that the designed networks drive the mean field approximation of the process to extinction quickly.  Since the seminal work of this approach was written \cite{preciado2014optimal}, authors have generalized the techniques to accommodate increasingly general single-process epidemic models \cite{nowzari2015optimal}, and multi-contagion epidemics \cite{watkins2015optimal}.  Unfortunately, it is not always the case that mean field approximation provides a good approximation to the stochastic epidemic process model; indeed, exponentially fast convergence of the mean-field model may not imply the same of its stochastic counterpart \cite{watkins2016optimal}.

Notably, literature considering approaches for explicitly incorporating feedback into the control of networked epidemic processes in the stochastic regime has been scarcely developed.  To the best knowledge of the authors, the only current approach is presented in the set of works \cite{drakopoulos2014efficient,drakopoulos2015lower,drakopoulos2015network}, which taken together present and analyze a method for controlling a restricted version of the continuous-time Susceptible-Infected-Susceptible (SIS) epidemic model.  In reference \cite{drakopoulos2014efficient}, the authors develop a suboptimal control policy for a continuous time $SIS$ process in which the healing rate of the process is varied as a function of the current state of the epidemic; however, computing the policy is computationally intractable.  In references \cite{drakopoulos2015lower} and \cite{drakopoulos2015network}, the authors provide analysis on the fundamental limitations of feedback control of continuous-time epidemic process, developing conditions on the budget available to the controller which are necessary and sufficient for the existence of a feedback control policy which causes extinction quickly.

Discrete-time networked epidemic processes have received considerably less attention than their continuous-time counterparts; the central themes in this area are captured by \cite{wang2003epidemic}, which establishes a stability condition for the mean field approximation of the discrete-time $SIS$ process, and \cite{wan2008designing} which uses heuristic design methods for designing static resource allocations which drive the epidemic to extinction.  More recent work has developed a method for designing controls via resource allocations when the underlying spreading network is not known, but a trace of the evolution of the mean field dynamics is available \cite{han2015data}.

\paragraph*{\textbf{Statement of Contributions}}
In this paper, we develop a novel method of feedback control for networked epidemic spreading processes, in particular the discrete-time $SIS$ process.  As the vast majority of prior work in the area of epidemic control is focused on open-loop, mean field control, our present work provides several considerable contributions. 
Chiefly, we believe this to be the first work to provide a computationally tractable method for introducing incomplete observations of the stochastic process in order to realize a feedback controller for networked epidemic processes.  Through the course of doing so, we provide the following contributions: (i) a novel model for observing epidemic processes and propagating the uncertainty in their dynamics (Section \ref{sec:observe}) and (ii) a computationally tractable method for feedback control of the $SIS$ process when a sufficiently rich subset of the nodes' compartmental memberships are observed (Section \ref{sec:RHC_SIS}).  Note that 
some proofs have been removed for purposes of brevity, and will be made available in separate, full-length publications.
\newgeometry{margin=0.75in} 
\subsection{Notation}
We denote by capital letters $X$ random variables.  When the random variable is a random vector, we let $X_i$ denote the $i$'th component of the random vector.   When the random variable is an element of a stochastic process $\{X(t)\},$ we will sometimes omit the time index, if it is clear from context.  Likewise, for some fixed time $t,$ we may reference the random variable $X(t)$ without a time index, the random variable $X(t+1)$ as $X^+,$ and the random variable $X(t-1)$ as $X^-.$  When $X$ is an indicator random variable, we will sometimes simply use $X$ to refer to the event $\{X = 1\},$ and the shorthand notation $X^c$ to refer to the event $\{X = 0\}.$  Accordingly, we use superscript $c$ to denote the complement.

\section{Modeling and Problem Statement}
We begin by reviewing basic concepts, definitions, and notations from the theory of the Susceptible-Infected-Suceptible ($SIS$) epidemic process. The material in this section is not novel, however it is necessary to place the contributions of this work properly in context, and so must be briefly developed to enable a self-contained discussion of the results.  

The $SIS$ process evolves on an $n$-node directed spreading graph, denoted by $\G = (V, E),$ where $V$ is the set of the graph's vertices, and $E$ is the ordered set of the graph's edges.  At each time $t$ in the process, each node $i$ in the graph takes on a value in one of two compartments, $S$ (for Susceptible), or $I$ (for Infected).  Transitions among the node's compartmental values occur stochastically; we denote the compartmental membership of node $i$ at time $t$ by a random variable $X_i(t),$ which takes the value $1$ when node $i$ is infected, and $0$ when node $i$ is susceptible.  Transitions from the infected compartment to the susceptible compartment for a particular node are not influenced by the compartmental memberships of other nodes in the graph, and occur with probability $\delta_i.$  Because it often reflects the probability of transitioning from sick to healthy, the value $\delta_i$ is commonly termed the healing rate of the process.  Transitions from the susceptible compartment to the infected compartment for a particular node are influenced by the compartmental memberships of the \emph{in-neighbors} of the node, and occur with probability $\prod_{j \in \Nin_i} (1- \beta_{ji}(t) X_j)$, where $\Nin_i$ denotes the set of in-neighbors of node $i$ in $\G,$ and the terms $\beta_{ji}$ are called the infection rates from node $j$ to node $i$.  Note that this value of infection probability reflects the fact that each infected neighbor $j$ attempts to infect node $i$ with probability $\beta_{ji}(t)$ at each time $t,$ with the infection taking root if and only if at least one of the infection attempts is successful.  All told, a stochastic description of the $SIS$ process is given by
\begin{equation}
	\label{eq:stochastic_SIS}
	X_{i}^+ = X_{i}(1-Y_{i}) + (1- X_{i})( 1 - Z_{i}(X)),
\end{equation}
in which $Y_{i}$ is an independent Bernoulli random variable with success probability $\delta_i,$ and $Z_{i}(X)$ is an independent Bernoulli random variable with success probability $\prod_{j \in \Nin_i} (1- \beta_{ji} X_j).$

As can be seen from \eqref{eq:stochastic_SIS}, the elements $\{X_i\}$ of the $SIS$ process exhibit complicated correlations: the future compartmental membership of a node depends explicitly on the compartmental membership of its in-neighbors, and so inexorably ties the evolution of the nodes' compartmental memberships together.  This is the central difficulty which has made control of the exact stochastic model of the process difficult.  Indeed, few works have considered control of the exact stochastic epidemic spreading problem, and those which have \cite{drakopoulos2014efficient,drakopoulos2015lower,drakopoulos2015network} consider performance analysis on controllers which assume that the controller has access to - and explicitly uses knowledge of - the compartmental membership of each node at every time in the process.  Our present work develops a method under which the epidemic need not be fully observed.

The analysis progresses in several distinct stages.  We propose and analyze an observation model of the $SIS$ process in Section \ref{sec:observe}, in which we show that if a sufficiently rich subset of vertices are observed, then the \emph{exact} conditional probabilities of the unobserved nodes can be propagated.  We build on this result in Section \ref{sec:RHC_SIS}, in which we develop a controller for the $SIS$ process which provably drives the process to extinction at a geometric rate in expectation.

\section{Observing the SIS Process}
\label{sec:observe}
In this section, we develop an observation model for the $SIS$ epidemic process~\eqref{eq:stochastic_SIS}. As our goal is to develop a tractable method which controls the process without approximating the epidemic dynamics, we can neither use the commonly studied mean field approximations, nor can we propagate the exact process dynamics forward for every node, as the former would introduce unacceptable approximation errors, and the latter would necessitate the use of $2^n$ equations to propagate the joint distribution \cite{hassibi2015SIRS}.
The approach we develop in this paper uses observations of the realized compartmental memberships of a subset of nodes as the process evolves in order to tractably make inferences and predictions.

\subsection{Concepts for Observation, Inference and Prediction}
We say a node $i \in V$ is observed if at each time $t,$ the compartmental membership of node $i$ is known to the controller.  We denote the set of nodes which the controller can observe by $\Obs \subseteq \V,$ and refer to it as the observer set.  It is intuitively clear that if a \emph{rich enough} subset of nodes are observed, then the compartmental membership probabilities of the unobserved nodes can be computed tractably.  To see that this must be the case, we trivially have tractable inferences if we assign every node to be observable, i.e. allowing $\Obs = \V.$  However, this is a na{\"i}ve approach: the entire process history would be observable.  In Section \ref{subsec:characterizing}, we establish structural properties of nontrivial observer sets which guarantee a desirable conditional independence property of the compartmental membership random variables.  To do so requires us to develop formal notation and definitions for observing, inferring, and predicting the $SIS$ process.

We denote by $\Hst_t$ the observed process history at time $t,$ i.e. the compartmental memberships $X_i(\tau)$ of all nodes $i \in \Obs$ for all times $0 \leq \tau \leq t.$  We denote by $\Theta_t$ the parameter history for the controlled $SIS$ process, i.e. the selected healing rates $\delta_i(\tau)$ and infection rates $\beta_{ij}(\tau)$ for all nodes and edges in the graph and all times $0 \leq \tau \leq t.$  
We consider the timing of the process to be such that at each time instant $\tau,$ the controller observes the compartmental memberships $\{X_i(\tau)\}_{i \in \Obs},$ updates its current estimates of the compartmental memberships of nodes not in the observer set, and then chooses a control action by selecting appropriate values for $(\vec{\delta}(\tau), \vec{\beta}(\tau)),$ i.e. the vectors containing the healing and infection rates for the process at time $\tau,$ based on predicted future compartmental memberships of the process. In order to enable a clear discussion of this process, we must define compact notation for making predictions about the future of the process and inferences of the process states.

Let $\Xhatitt$ denote the conditional probability that node $i$ is infected at time $t$, given the information available to the controller at time $t$ prior to the selection of the parameters at time $t,$ i.e. $\Xhatitt \triangleq \Pr(X_i(t) = 1 | \Hst_t, \Theta_{t-1}).$  This represents our current belief about the compartmental membership of node $i,$ given all available information to the controller at time $t.$  Similarly, let $\Xhatitpot$ denote the conditional probability that node $i$ is infected at time $t+1$ given the information available to the controller at time $t,$ including the design of the parameters at time $t,$ i.e. $\Pr(X_i(t+1) = 1 | \Hst_t, \Theta_{t})$.  This is a prediction of the next compartmental membership of node $i,$ given the history through time $t.$

Understanding which observer sets enable the tractable computation of $\Xhatitt$ and $\Xhatitpot$ is important, and nontrivial.  To understand the central issue which makes this computation difficult, note that the stochastic $SIS$ process evolves in the $2^n$ dimensional space of all combinations of compartmental memberships, and so using standard inference algorithms for graphical models (see, e.g. \cite{koller2009probabilistic}) to estimate the joint distributions of the process would not be tractable.  Moreover, marginalizing over the joint distribution to recover the estimate $\Xhatitt$ would require a sum over the set of $O(2^n)$ dimension which contains all possible states of the process.  We address this issue in the following subsection.

\subsection{Characterizing Efficient Observer Sets}
\label{subsec:characterizing}
In this section, we use tools from graphical models to provide rigorous conditions under which the compartmental membership random variables $\{X_i(t)\}$ are conditionally independent, given the information the controller has through time $t.$  This property is of central importance to enabling tractable prediction and inference of the controlled $SIS$ process, and so must be developed fully in order to continue properly.  To state the result concisely, we first need to define the concept of the moralized spreading graph.  For the reader familiar with the theory of graphical models in statistical inference, it will be recognized as a direct extension of the moralization of a directed graphical model.
\vspace{-4pt}
\begin{defn}[Moralized Spreading Graph]
	Consider a spreading graph $G.$  The moralization of $G$ is defined as the undirected graph $\mor(G) = (V,\Emor),$ where 
	\begin{equation*}
	\begin{aligned}
	\Emor \triangleq &\{\{i,j\} \in V \times V \, | \, (i,j) \in E\} \, \cup \\
	&\{\{i,k\} \in V \times V \, | \, (i,j) \in E, (k,j) \in E\}.
	\end{aligned}  
	\end{equation*}
\end{defn}
Informally, $\mor(G)$ is the graph which results by connecting any two nodes which share a common out-neighbor, and removing directionality from the graph.  Intuitively, the edges of $\mor(G)$ are structured so as to represent the structure of relationships among the random variables of the process.

We now state our first major result, which demonstrates the importance of the moralized spreading graph in enabling tractable feedback in networked epidemic control.  In particular we show that taking an observer set $\Obs$ which forms a \emph{vertex cover} of $\mor(G)$, i.e. a collection $\C$ of nodes such that each edge $e = \{v,v^\prime \}$ of $\mor(G)$ has at least one of $v$ or $v^\prime$ in $\C,$ is necessary and sufficient for a desirable conditional independence property to hold.  When this property holds, we can perform inference and prediction tractably, and hence our observer set is sufficiently well-structured.
\vspace{-8pt}
\begin{theorem}[Conditional Independence of $SIS$]
	\label{thm:cond_ind}
	Consider an instance of the $SIS$ spreading process on $G = (V,E).$  The random variables $\{X_i(t)\},$ are mutually independent conditioned on the process history $\Hst_t$ and the parameter history $\Theta_{t-1}$ for every sequence of spreading parameters if and only if the observer set $\Obs \subseteq V$ forms a vertex cover of the moralized spreading graph, $\mor(G).$ 
\end{theorem}
\vspace{-3pt}

We close this section by informally making some observations on the size of observer sets which satisfy Theorem \ref{thm:cond_ind}.  For the $k$-node star graph with directed edges from the hub to the leaves, only the hub needs to be observed.  For the $k$-node complete graph, any collection of $k-1$ nodes will suffice as an observer set, but all vertex covers must contain at least $k-1$ nodes.  Hence, the fraction of observed nodes may tend to zero or to one as the number of nodes in the graph tends to infinity, depending on the topology of the graph.  These examples demonstrate that high directionality in a graph allows for simple observability, but high density of edges makes observation difficult. 
  
Note also that in general, computing minimal cardinality observer sets which satisfy Theorem \ref{thm:cond_ind} is a difficult problem.  In fact, it is equivalent to the vertex cover problem.  We will not go into detail here, but this equivalence allows us to prove that the problem is $\NP$-hard, but easy to approximate \cite{williamson2011design}.
\vspace{-5pt}
\subsection{Inference and Prediction of the Conditional SIS Process}
\label{subsec:inference}

In this section, we develop tractable methods for inference and prediction under the assumption that the conditions of Theorem \ref{thm:cond_ind} are satisfied.  We first state the result formally, then use the remainder of the section for proof and discussion.
\vspace{-7pt}
\begin{theorem}[Inference and Prediction]
	\label{thm:condSIS}
	Consider an instance of the $SIS$ process on $G.$  If the observer set chosen forms a vertex cover of the moralized spreading graph, then the conditional inference probabilities $\Xhatitt,$ and prediction probabilities $\Xhatitpot$ can be computed using only the state estimates and spreading parameters of time $t-1,$ and the observations made and parameters selected at time $t.$  
	
	More precisely, the worst-case complexity of computing each term is $O(\dmax^2),$ where $\dmax$ is the maximum in-degree of the spreading graph.
\end{theorem} 
We now prove Theorem \ref{thm:condSIS} by constructing state estimators and predictors with the specified properties.  For the convenience of the reader, we do so systematically, with each constructed mechanism appearing in its own subsection.

\subsubsection*{\textbf{Inference of an Observed Node}}
It is intuitive that in this case, $\Xhatitt = X_i(t)$ for all $t,$ and so the update can be performed in constant time and only requires the observation of the process at time $t.$  Formally, this follows from noting that for all $t,$ $X_i(t)$ is measurable with respect to $\I_t = \{\Hst_t, \Theta_{t} \},$ and so is known to the controller. 

\subsubsection*{\textbf{Inference of an Unobserved Node}}
This is the most technically challenging case, and requires the most involved calculation.  Fundamentally, our task is in showing that the standard Bayes' rule expression
\begin{equation*}
	\Pr(X_i| \I_{t-1}, \{X_k\}_{k \in \Obs}) = \frac{\Pr(X_i,\{X_k\}_{k \in \Obs} | \I_{t-1})}{\Pr(\{X_k\}_{k \in \Obs} | \I_{t-1})}, \, i \notin \Obs
\end{equation*}
can be evaluated in a computationally tractable manner, and by only using the state estimates and parameters from time $t-1$ together with the observations made at time $t,$ as claimed by the statement of the theorem.

By considering Theorem \ref{thm:cond_ind}, it may be shown that the collection of random variables which are not the compartmental memberships of in-neighbors of $i$, i.e. $\{X_j\}_{j \in \Obs \cap V\setminus\Nin_i},$ is conditionally independent of $X_i,$ given the compartmental memberships of in-neighbors, i.e. $\{X_j\}_{j \in \Obs \cap \Nin_i},$ and the process history.  It then follows that when updating the estimate of node $i,$ only observations of $i$'s in-neighbors effect the estimate; more formally, 
\begin{equation}
\label{eq:Obs_comp_step1}
\begin{aligned}
\Pr(X_i|\{X_j\}_{j \in \Obs}, \I_{t-1}) = \Pr(X_i|\{X_j\}_{j \in \Obs \cap \Nin_i}, \I_{t-1})
\end{aligned}
\end{equation}
holds.  We may refine \eqref{eq:Obs_comp_step1} further by noting that the same argument holds when applied to in-neighbors of $i$ which were observed to have transitioned from an infectious state to a healthy state.  In particular, the transitions of in-neighbors of node $i$ which are infected at time $t-1$ are determined independently from all other transitions, and so the identity
\begin{equation}
\label{eq:xhatitt_obs_comp1}
\begin{aligned}
\Xhatitt = \frac{\Pr(X_i, \{X_k\}_{\Kihx}| \I_{t-1})}{\Pr(\{X_k\}_{\Kihx} | \I_{t-1})},
\end{aligned}
\end{equation}
holds, where we have defined $\Kihx$ to be the index set of the in-neighbors of $i$ which were observed to have been healthy at time $t-1.$  It remains to show that the conditional probabilities $\Pr(X_i, \{X_k\}_{\Kihx}| \I_{t-1})$ and $\Pr(\{X_k\}_{\Kihx} | \I_{t-1})$ can be computed tractably, and by only using the state estimates and parameter values at time $t-1.$  

An application of the law of total probability gives
\begin{equation}
\label{eq:xhatitt_obs_comp2}
\begin{aligned}
\Pr(&X_{i,}\{X_k\}_{\Kihx} | \I_{t-1}) = \\
&\Pr(X_{i,}\{X_k\}_{\Kihx} | \I_{t-1}, X_i^{-}) \Xhatitmotmo \, +\\
&\Pr(X_{i,}\{X_k\}_{\Kihx} | \I_{t-1}, X_i^{c-}) (1 - \Xhatitmotmo),
\end{aligned}
\end{equation}
\begin{equation}
\label{eq:xhatitt_obs_comp3}
\begin{aligned}
\Pr(&\{X_k\}_{\Kihx} | \I_{t-1})  = \\
&\Pr(\{X_k\}_{\Kihx} | \I_{t-1}, X_i^{-}) \Xhatitmotmo\, +\\
&\Pr(\{X_k\}_{\Kihx} | \I_{t-1}, X_i^{c-}) (1 - \Xhatitmotmo),
\end{aligned}
\end{equation}
and the conditional independence of healing and infection events affecting node $i$ at time $t-1,$ given $\I_{t-1}$ and the compartmental membership of node $i$ at time $t-1$ gives
\begin{equation}
\label{eq:xhatitt_obs_comp4}
\begin{aligned}
&\Pr(X_{i,}\{X_k\}_{\Kihx} | \I_{t-1}, X_i^{-}) = \\
&\;\;\;\; (1 - \delta_i^{-})\Pr(\{X_i\}_{\Kihx} | \I_{t-1}, X_i^{-}),\\
&\Pr(X_{i,}\{X_k\}_{\Kihx} | \I_{t-1}, X_i^{c-}) = \\
&\;\;\;\; (1 - \Pi_{j \in \Nin_i} (1 - \beta_{ji}^{-} X_j^{-})) \Pr(\{X_i\}_{\Kihx} | \I_{t-1}, X_i^{c-}).
\end{aligned}
\end{equation}
Noting that the decomposition given by \eqref{eq:xhatitt_obs_comp4} has still not removed the effect of the process history, we must complete the computation by 
defining the notation $\Kihh$ as the in-neighbors of node $i$ which were healthy at time $t-1$ and are healthy again at time $t,$ and $\Kihi$ as the in-neighbors of node $i$ which were healthy at time $t-1,$ and are infected at time $t,$ and noting
\begin{equation}
\label{eq:xhatitt_obs_comp5}
\begin{aligned}
&\Pr(\{X_k\}_{\Kihx} | \I_{t-1}, X_i^{-}) = \\
&\;\;\;\;\Pi_{k \in \Kihh} \left((1 - \beta_{ik}^-) \Pi_{j \in \Nin_k \setminus \{i\}}(1 - \beta_{jk}^- X_j^-) \right)\\
&\;\;\;\;\Pi_{k \in \Kihi} \left(1 - (1 - \beta_{ik}^-)\Pi_{j \in \Nin_k \setminus \{i\}}(1 - \beta_{jk}^{-} X_j^{-}) \right),\\
&\Pr(\{X_k\}_{\Kihx} | \I_{t-1}, X_i^{c-}) = \\
&\;\;\;\;\Pi_{k \in \Kihh} \left(\Pi_{j \in \Nin_k \setminus \{i\}}(1 - \beta_{jk}^- X_j^-) \right)\\
&\;\;\;\;\Pi_{k \in \Kihi} \left(1 - \Pi_{j \in \Nin_k \setminus \{i\}}(1 - \beta_{jk}^{-} X_j^{-}) \right),
\end{aligned}
\end{equation}
where we have implicitly used the facts that if the observer set forms a vertex cover of the moralized spreading graph, then an unobserved node's in-neighbors are necessarily observed, and any observed node $k$ has at most one unobserved in-neighbor $i$ to explicitly remove the effect of the process history on the computed estimate.  Note that in an implementation, equations \eqref{eq:xhatitt_obs_comp1}-\eqref{eq:xhatitt_obs_comp5} must be computed in \emph{reverse} order, but the complexity of the computation is at worst $O(\dmax^2),$ and so verifies the theorem's claim.

\subsubsection*{\textbf{Prediction of an Observed Node}}
The computation in this case is straightforward.  It can be shown that if $i$ is in $\Obs,$ it must be that at most one in-neighbor of $i$ is not in $\Obs.$  Indeed, if two nodes $j$ and $k$ are neighbors of $i$ in $\mor(G),$ they are also neighbors of each other in $\mor(G),$ or at least one of $j$ or $k$ is not an in-neighbor of $i$ in $G.$

We may use this fact to factor the calculation as
\begin{equation}
\label{eq:xhatitpot_obs}
\begin{aligned}
&\Xhatitpot = X_i(1- \delta_i)+ (1- X_i) &\\
&\hspace{10pt}\left(1 - (1- \beta_{j^\prime i} \hat{x}_{j^\prime}(t | t)) \Pi_{j \in \Nin_i \cap \Obs} (1 - \beta_{ji} X_j)\right),
\end{aligned}
\end{equation}
where we define $j^\prime$ as the unique unobserved in-neighbor of $i$ in $G.$  Noting that \eqref{eq:xhatitpot_obs} requires $O(d_{\max})$ operations and is a function of only the current estimates and selected parameters verifies that the predictor has the claimed properties.

\subsubsection*{\textbf{Prediction of an Unobserved Node}}
The computation in this case is not much different than that of the predictions made for observed nodes.  Indeed, by conditioning on the compartmental membership of node $i$ at time $t$ and evaluating the remaining conditional probabilities, we arrive at
\begin{equation}
\label{eq:xhatitpot_obs_comp}
\begin{aligned}
&\Xhatitpot = (1 - \delta_i) \Xhatitt+&\\
&\hspace{55pt}(1 - \Pi_{j \in \Nin_i} (1 - \beta_{ji} X_j))(1 - \Xhatitt),
\end{aligned}
\end{equation}
where we have implicitly used the fact that if $\Obs$ forms a vertex cover of the moralized spreading graph, it implies that the in-neighbors of unobserved nodes are observed.  Here again, the complexity of the computation required is $O(\dmax),$ and  its evalutation only requires the state estimates and selected parameters at time $t,$ and so satisfies the claimed properties.  Moreover, this realization completes our proof. 

It is important to note that the designed estimation and prediction algorithms produced unbiased estimates and predictions, provided they begin with an unbiased estimate of the initial compartmental memberships.  This can be proven formally by noting that the above computations carried the relevant conditional expectations through exactly, and so remain unbiased if they began as such. 
\vspace{-1pt}
\section{Rolling Horizon Control of Conditional SIS}
	\label{sec:RHC_SIS}
\vspace{-1pt}
	In this section, we consider the problem of minimizing the cost realized by a controller tasked with driving the $SIS$ process so as to maintain a minimum decay rate in expectation of the total infection over the entire population.  Let $f_i$ be the cost function associated to the healing rate of node $i,$ $g_{ij}$ be the cost function associated to the infection rate of the edge $(i,j),$ and $r$ be a decay constant in the open unit interval.  Then, we may write the problem we study in this section as the chance constrained optimization problem:
\begin{equation}
	\label{prob:RHC_SIS}
	\begin{aligned}
	& \underset{\vec{\delta},\vec{\beta}}{\text{minimize}}
	& & \sum_{i \in V} f_i(\delta_i(t))+\sum_{(i,j) \in E} g_{ij}(\beta_{ij}(t)) \\
	& \text{s.t.}
	& & \sum_{i \in V} \bE \left[X_i(t+1) | \I_t \right] \leq r \sum_{i \in V} \bE \left[X_i(t) | \I_t \right],
	\end{aligned}
\end{equation}
in which the objective function is the cost incurred by the parameters chosen by the controller at time $t,$ and the constraint is enforced to ensure that the closed-loop system is such that the number of infected nodes decays in expectation with at least a rate $r$ at every time step.  We begin our analysis of this problem by showing that under mild assumptions, a coordinate change can be utilized to pose \eqref{prob:RHC_SIS} as an equivalent convex program.  We now state this formally.
\vspace{-3pt}
\begin{theorem}
	[Convex Bayesian SIS Control]
	\label{thm:RH_convex}
	Consider the convex optimization program
	\begin{equation}
	\label{prob:RHC_SIS_convex}
	\begin{aligned}
	& \underset{\vec{\delta}^c,\vec{\gamma}}{\text{minimize}}
	& & \sum_{i \in V} f_i(\delta_i^c)+\sum_{(i,j) \in E} \tilde{g}_{ij}(\gamma_{ij}) \\
	& \text{s.t.}
	& & \sum_{i \in \Obs} \delta_i^c X_i + \psiiw(X,\xhat, \vec{\gamma}) X_i^c +\\
	&&&\sum_{i \in \Obs^c} \delta_i^c \Xhatitt + \psiiw(X,\xhat, \vec{\gamma}) \hat{x}^c_i(t|t)\\
	&&&\leq r \sum_{i \in V} \Xhatitt,\\
	\end{aligned}
	\end{equation}
	where we additionally restrict the variables $\delta_i^c$ and $\gamma_{ij}$ to the closed unit interval, and have defined the convex functions
	\begin{equation}
	\label{eq:SIS_convex}
	\begin{aligned}
	&\psiiw(X,\xhat,\vec{\gamma})=&
	&\begin{cases}
	1-  \hat{x}_{j^\prime}(t|t) \gamma_{j^\prime i}^{\frac{1}{w}} \Pi_{j \in \{\Nin_i \cap \mathcal{X}_i\}} \; \gamma_{ji}^{\frac{1}{w}}\\
	\; \; \;  -\hat{x}_{j^\prime}^c(t | t) 
	\Pi_{j \in \{\Nin_i \cap \mathcal{X}_i \}} \gamma_{ji}^{\frac{1}{w}},&\, i \in \Obs\\
	1 - \Pi_{j \in \{\Nin_i \cap \mathcal{X}_i \}} \;  \gamma_{ji}^{\frac{1}{w}},&\, i \in \Obs^c
	\end{cases}
	\end{aligned}
	\end{equation}
	where $w > \dmax,$ the sets $\mathcal{X}_i = \{i \in V \cap \Obs \, | \, X_i = 1\}$, and the shorthand notation $X_i^c = (1 - X_i),$ $\hat{x}_i^c(t | t) = (1 - \hat{x}_i(t|t)),$ and $\delta_i^c = (1 - \delta_i)$ for purposes of compacting notation.  Suppose the functions $f_i$ and $\tilde{g}_{ij} = g_{ij}(1 - \gamma_{ij}^{\frac{1}{w}} )$ are convex in the variables $\delta_i^c,$ and $\gamma_{ij},$ respectively.  Then \eqref{prob:RHC_SIS} and \eqref{prob:RHC_SIS_convex} are equivalent optimization problems, where the optimal healing rates of \eqref{prob:RHC_SIS} can be computed as $\delta_i^\star = 1 - \delta_i^{c \star},$ and the optimal infection rates of \eqref{prob:RHC_SIS} can be computed as $\beta_{ij} = 1 - (\gamma_{ij}^\star)^{\frac{1}{w}},$ where $\delta_i^{c \star}$ and $\gamma_{ij}^\star$ are solutions to \eqref{prob:RHC_SIS_convex}.
\end{theorem}

\begin{remark}[Convexity of $\psiiw$ Functions]
	{\rm Note that the convexity of each $\psiiw$ can be verified by applying an established result from signomial optimization literature \cite{lundell2009transformation}.  However, it is worth noting that product terms are in general nonconvex, and so CVX \cite{CVX} may not solve the problem.  Solutions can be obtained by coding standard convex optimization algorithms (see, e.g. \cite{nocedal2006numerical}).} \oprocend
\end{remark}

\begin{remark}[Convexity of Objective Functions]
	{\rm
It is worth noting that the class of objective functions which are convex under the proposed variable transformation may not be as large as desired.  While this is a technical restriction, the importance of Theorem \ref{thm:RH_convex} remains.  The transformation convexifies the region defined by the decay rate constraint, and so finding a feasible suboptimal controller is always tractable, irrespective of the objective function, and so suboptimal control is always computationally tractable. \oprocend
}
\end{remark}
The proposed controller causes the controlled $SIS$ process to converge to the all-healthy state quickly in expectation; the formal statement of our convergence result follows.
\begin{theorem}
	[Convergence of Controlled SIS]
	\label{thm:RH_convergence}
	For any $r \in (0,1)$, the controlled $SIS$ process converges to the disease-free equilibrium in expectation geometrically at rate $r.$
\end{theorem} 

We close this section by noting that the performance results presented in this paper can be refined to demonstrate stronger properties, such as an upper bound on the expected time-to-extinction for the controlled process, and almost sure finite time convergence.  However, the analysis required to formalize these features is too lengthy to discuss in detail in this venue, and so is saved for publication in future work.
\vspace{-5pt}
\section{Simulations}
\label{sec:sims}
\vspace{-5pt}
In this section, we use extensive numerical simulation in order to verify the utility of our analytical findings.  For purposes of simplicity, we chose our objective functions to be $f_i(\delta_i) = \delta_i,$ and $g_{ij}(\beta_{ij}) = (1 - \beta_{ij})^{\dmax-1},$ which under transformation becomes $g_{ij}(\gamma_{ij}) = \gamma_{ij}^{\frac{\dmax-1}{w}}.$  In this paper, we report the results of a study of a random graph with $30$ nodes, connection probability $p = 0.2,$ max degree $11,$  and an observer set of $24$ nodes.

The results are given in Figures \ref{fig:x_hist} and \ref{fig:obj_val_plot}.  Inspection of Figure \ref{fig:x_hist} demonstrates that convergence of the controlled $SIS$ process occurs at a geometric rate, as demonstrated in the analysis.  Inspection of Figure \ref{fig:obj_val_plot} reveals a structure that experience has revealed to be typical of the realized controllers.  When the graph is in a highly infectious or a nearly healthy state, the control cost is low.  When the graph is in a middling state, the cost is high.  This is explained by noticing that when a graph is in a highly infectious state, only the healing rates need to be controlled to guarantee a particular decay rate.  Likewise, when a graph is nearly healthy, only infection rates need to be aggressively controlled to guarantee decay.  In between these extremes, both sets of resources need to be utilized, resulting in a higher cost.  

{	
	\begin{figure}[htb] 
		\begin{center}
			{\includegraphics[width = 0.45\textwidth]{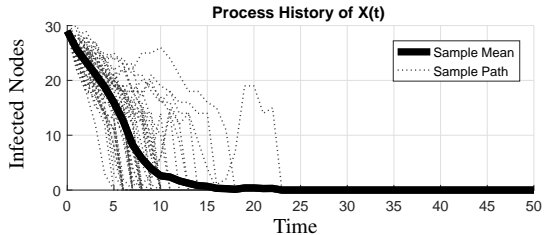}}
			\put(-120,-7){\footnotesize Time}
			\put(-220,17){\footnotesize \rotatebox{90}{Infected Nodes} }
		\end{center}
		\caption{A study of the convergence of the controlled $SIS$ process.  Time is plotted on the horizontal axis, and the number of infected nodes on the vertical axis, with dotted lines representing individual sample paths, and the dark line the sample mean.} \label{fig:x_hist}
	\end{figure}
}

{	
	\begin{figure}[htb] 
		\begin{center}
			{\includegraphics[width = 0.45\textwidth]{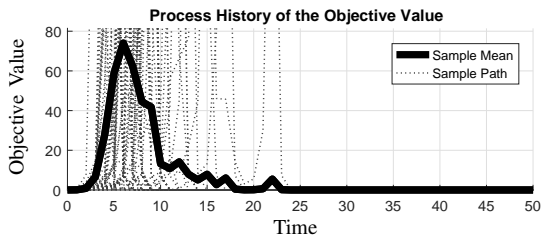}}
			\put(-120,-7){\footnotesize Time}
			\put(-220,15){\footnotesize \rotatebox{90}{Objective Value}}
		\end{center}
		\caption{A study of the cost of allocated resources to control the $SIS$ process.  Time is plotted on the horizontal axis, and the realized value of the objective function is plotted on the vertical axis, with dotted lines representing individual sample paths and the dark line representing the sample mean of the simulated trajectories.} \label{fig:obj_val_plot}
	\end{figure}
} 

The dynamic nature of the allocated budget is a significant advantage over mean field approaches.  In the limit where the estimated state of the process approaches the zero infection state, the cost of the realized control approaches zero as well.  This suggests that for an epidemic which extincts in finite time, only a finite budget is needed, which is indeed infinitely better than the static controllers realized by mean field methods.  Note that in a mean field approach, the allocated resources remain static throughout the evolution of the process, which implies that the integrated cost of the controller diverges as time progresses.  Introducing feedback into the process avoids expending resources needlessly.  

\section{Summary and Future Work}
	
There are several possible directions for future work, but a few which seem of particular importance.  In an application of networked epidemic theory - be it biological, technological, or theoretical - it seems likely that the systems under consideration will be very large, and subjected to uncertainties beyond the scope of those considered here.  As such, it would be worthwhile to investigate methods for incorporating stochastic healing and infection rates.  Moreover, it may be the case that the observability conditions outlined in this text cannot be attained, and so it may be of interest to consider means for tractably incorporating information of the process, even when inference and prediction may only be done approximately.  It seems possible that such advances can be made with continued effort by the research community.
\section*{Acknowledgments}

This is supported by the TerraSwarm Research Center, one of six centers supported by the STARnet phase of the Focus Center Research Program (FCRP), a Semiconductor Research Corporation program sponsored by MARCO and DARPA.

\bibliography{research}
\bibliographystyle{ieeetr}
\end{document}